\documentclass[12pt,leqno]{amsart}
\usepackage{amssymb}
\usepackage{amsmath,amssymb,color}
\numberwithin{equation}{section}

\newcommand{\eexx}{\vert \exp(\zeta^{\frac{1}{\alpha}}) \vert}
\newcommand{\ep}{\varepsilon}
\newcommand{\la}{\lambda}
\newcommand{\va}{\varphi}
\newcommand{\ppp}{\partial}

\newcommand{\pppa}{\partial_t^{\alpha}}

\newcommand{\sumij}{\sum_{i,j=1}^d}
\newcommand{\uab}{u_{a,b}}

\newcommand{\R}{\mathbb{R}}

\newcommand{\C}{\mathbb{C}} 
\newcommand{\N}{\mathbb{N}}

\newcommand{\ooo}{\overline}
\newcommand{\OOO}{\Omega}
\newcommand{\MLO}{E_{\alpha,1}}
\newcommand{\MLT}{E_{\alpha,2}}
\newcommand{\MLAA}{E_{\alpha,\alpha}}
\newcommand{\suml}{\sum_{j=1}^{\ell_n}}
\newcommand{\sumn}{\sum_{n=1}^{\infty}}

\allowdisplaybreaks

\title
[]
{
Backward problems in time for fractional diffusion-wave equation}

\author{
$^1$ G. Floridia
and \, $^{2,3,4}$ M.~Yamamoto }

\thanks{
$^1$ 
Department PAU, 
Universit\`a Mediterranea di Reggio Calabria,
Via dell'Universit\`a 25  
89124 Reggio Calabria, Italy \&
INdAM Unit, University of Catania, Italy,
{\tt floridia.giuseppe@icloud.com} \\
$^2$ Graduate School of Mathematical Sciences, The University
of Tokyo, Komaba, Meguro, Tokyo 153-8914, Japan \\
$^3$ Honorary Member of Academy of Romanian Scientists, 
Splaiul Independentei Street, no 54,
050094 Bucharest Romania \\
$^4$ Peoples' Friendship University of Russia 
(RUDN University) 6 Miklukho-Maklaya St, Moscow, 117198, Russian Federation
e-mail: {\tt myama@ms.u-tokyo.ac.jp}
}

\date{}
\begin{document}
\maketitle

\baselineskip 18pt

\begin{abstract}
In this article, for a time-fractional diffusion-wave equation 
$\pppa u(x,t) = -Au(x,t)$, $0<t<T$ with fractional order $\alpha \in (1,2)$, 
we consider the backward problem in time:
determine $u(\cdot,t)$, $0<t<T$ by $u(\cdot,T)$ and 
$\ppp_tu(\cdot,T)$.  We proved that there exists 
a countably infinite set $\Lambda \in (0,\infty)$ with  
a unique accumulation point $0$ such that the backward problem 
is well-posed for $T \not\in \Lambda$.
\end{abstract}

\section{Introduction and main results}

Let $\OOO$ be a bounded domain in $\R^d$ with sufficiently smooth 
boundary $\ppp\OOO$.  We consider a fractional differential equation:
$$
\pppa u(x,t) = -\mathcal{A}u(x,t), \quad x \in \OOO, \,  0<t<T, 
$$
where $-\mathcal{A}$ is a uniformly elliptic operator.  Henceforth for
$n-1 < \alpha < n$ with $n\in \N$, we define the Caputo derivative by 
$$
\pppa g(t) = \frac{1}{\Gamma(n-\alpha)} \int^t_0 (t-s)^{n-\alpha-1}
\frac{d^n}{ds^n}g(s) ds.
$$
For $\alpha \in (0,1) \cup (1,2)$, equation (1.1) is widely studied
not only by mathematical interests but also for the modelling of 
various types of diffusion phenomena in heterogeneous media.
Among them, we particularly refer to the anomalous diffusion which cannot
be modelled by a classical advection-diffusion equation which corresponds to 
$\alpha=1$.  More precisely, field data of diffusion of e.g., contaminants 
in soil often indicate long-tailed profiles, which cannot be interpreted 
by a classical advection-diffusion equation whose solution decays very fast,
i.e., exponentially.

There are tremendously many works on mathematical analysis and here we 
are strongly limited to some references.  As for the well-posedness of 
the initial boundary value problem for (1.1), we refer to 
Kubica, Ryszewska and Yamamoto \cite{KRY}, Kubica and Yamamoto \cite{KY},
Sakamoto and Yamamoto \cite{SY}, Zacher \cite{Za}, and for inverse problems
and related topics the readers can consult the handbook Li, Liu and 
Yamamoto \cite{LLY}, Li and Yamamoto \cite{LiY}, Liu, Li and 
Yamamoto \cite{LiuLiY}.

A solution to equation (1.1) with $\alpha \ne 1$
shows behavior which is 
essentially different from the case of $\alpha = 1$
and can characterize the anomaly of the diffusion in the heterogeneous
media.  Among such characteristic properties, the backward stability in time
is important and this is the main subject of this article.
In the case of $\alpha=1$, the classical diffusion equation possesses 
the strong smoothing property, so that we cannot solve the equation 
with final value condition, and cannot have good stability 
but with given a priori bound assumptions,
one can prove only conditional stability of logarithmic type 
(e.g., Imanuvilov and Yamamoto \cite{IY}, Isakov \cite{Is}, 
Section 9 in Yamamoto \cite{Ya}).
 
In the case of $0<\alpha<1$, Sakamoto and Yamamoto
\cite{SY} established the well-posedness of the backward problem in time
under reasonable regularity condition.  After \cite{SY}, as for 
$0 < \alpha < 1$, there have been many theoretical and numerical works 
on the backward problems, and here we can refer to 
Floridia, Li and Yamamoto \cite{FLY}, 
Liu and Yamamoto \cite{LY}, 
Tuan, Huynh, Ngoc and Zhou \cite{THNZ}, Tuan, Lung and Tatar \cite{TLT},
Tuan, Thach, O'Regan and Can \cite{TNOZ}, 
Wang, Wei and Zhou \cite{WWZ}, Wang and Liu \cite{WL2},
Wei and Wang \cite{WW}, 
Xiong, Wang and Li \cite{XWL}, Yang and Liu \cite{YL}, and we do not 
intend comprehensive references.

However, to the best knowledge of the authors, except for
Wei and Zhang \cite{WeiZ}, there are still no works on the 
backward problem as long as the case $1<\alpha<2$ is concerned, although
the case $1<\alpha<2$ is used for the modelling.

The purpose of this article is to sharpen the stability and the 
uniqueness, which improves the theoretical achievements of \cite{WeiZ}
for the backward problem for the case of $1<\alpha<2$.
\\

For the formulation of the problem, we introduce an operator
and function spaces.
We assume that all functions under consideration are real-valued.
Henceforth $L^2(\OOO)$ and $H^2(\OOO)$, $H^1_0(\OOO)$, $H^2_0(\OOO)$, etc.
denote the Lebesgue space and usual Sobolev spaces (e.g.,
Adams \cite{Ad}), and by $\Vert \cdot\Vert_X$
we denote the norm in the space $X$.  We set $(a,b) = \int_{\OOO} a(x)b(x) dx$.
Identifying the dual space $(L^2(\OOO))'$ with 
itself, we denote $H^{-1}(\OOO) = (H^1_0(\OOO))'$ and
$H^{-2}(\OOO) = (H^2_0(\OOO))'$.

We set
$$
-\mathcal{A}v(x) = \sumij \ppp_i(a_{ij}(x)\ppp_jv(x)) 
+ c(x)v(x), \quad x \in \OOO
$$
for $v \in C^2(\ooo{\OOO})$, where $\ppp_i = \frac{\ppp}{\ppp x\i}$ for
$1\le i \le d$, 
$a_{ij} = a_{ji} \in C^1(\ooo{\OOO})$, $c\in C(\ooo{\OOO})$, $c\le 0$ on 
$\ooo{\OOO}$.
Then we define an operator $A$ in $L^2(\OOO)$ by
$$
Av = \mathcal{A}v, \quad v \in \mathcal{D}(A) 
:= H^2(\OOO) \cap H^1_0(\OOO).
                                       \eqno{(1.2)}
$$
Here $v \in \mathcal{D}(A)$ means that $v=0$ on $\ppp\OOO$ 
in the sense of the trace.
Then it is also known that the operator defined by (1.2) 
has eigenvalues and we number the set of 
all the eigenvaleus:
$$
0 < \mu_1 < \mu_2 < \cdots \longrightarrow \infty.
$$
Let $\{ \va_{nj}\}_{1\le j \le \ell_n}$ be an orthonormal basis of
Ker $(A-\mu_n)$: $A\va_{nj} = \mu_n\va_{nj}$ and $(\va_{nj}, \va_{mi})
= \delta_{nm}\delta_{ij}$ where we set $\delta_{ij} = 1$ if $i=j$ and $=0$ if 
$i \ne j$.  Then we see that $\{ \va_{nj};\, n\in \N, \, 1 \le j \le \ell_n
\}$ is an orthonormal basis in $L^2(\OOO)$.

Throughout this article, we always assume 
$$
1 < \alpha < 2.
$$

In terms of $A$, we rewrite (1.1) as 
$$
\left\{ \begin{array}{rl}
& \pppa u(x,t) = -Au(x,t), \quad x\in \OOO, \,
t>0, \\
& u(\cdot,0) = a, \quad \ppp_tu(\cdot,0) = b, \qquad x \in \OOO. \\
\end{array}\right.
                                          \eqno{(1.3)}
$$

By $E_{\alpha,\beta}(z)$ we denote the Mittag-Leffler function:
$$
E_{\alpha,\beta}(z) = \sum_{k=0}^{\infty} \frac{z^k}{\Gamma(\alpha k + \beta)},
$$
with $\alpha>0$ and $\beta \in \C$, $z\in \C$.  It is known that
$E_{\alpha,\beta}(z)$ is an enire function in $z\in \C$ (e.g.,
Gorenflo, Kilbas, Mainardi and Rogosin \cite{GKMR}, Podlubny \cite{Po}).
 
Before stating the main results, we show the well-posedness and 
the regularity of the solution $u_{a,b}$ to (1.3).  
\\
\vspace{0.1cm}
{\bf Proposition.}\\
{\it 
Let $a, b\in L^2(\OOO)$.  Then there exists a unique 
solution $\uab$ to (1.3) such that
$$
\left\{ \begin{array}{rl}
& \uab  \in C([0,T];L^2(\OOO)) \cap C((0,T;H^2(\OOO) \cap H^1_0(\OOO)),\\
& \lim_{t\to 0} \Vert u(\cdot,t) - a\Vert_{L^2(\OOO)} 
= \lim_{t \to 0} \Vert \ppp_tu(\cdot,t) - v\Vert_{H^{-2}(\OOO)} = 0
\end{array}\right.
$$
and 
$$
\left\{ \begin{array}{rl}
& u(x,t) = \sumn \suml \{ (a,\va_{nj})\MLO(-\mu_nt^{\alpha})
+ (b,\va_{nj}) t\MLT(-\mu_nt^{\alpha})\} \va_{nj}(x)\\
&\qquad \qquad \qquad \qquad \qquad \qquad \qquad \qquad \qquad
\qquad \qquad \qquad \mbox{in $C([0,T];L^2(\OOO))$}, \\
& \ppp_tu(x,t) = \sumn \suml \{ 
-\mu_nt^{\alpha-1}(a,\va_{nj})\MLAA(-\mu_nt^{\alpha})
+ (b,\va_{nj}) \MLO(-\mu_nt^{\alpha})\} \va_{nj}(x) \\
&\qquad \qquad \qquad \qquad \qquad \qquad \qquad \qquad \qquad
\qquad \qquad \qquad \mbox{in $C([0,T];L^2(\OOO))$}.
\end{array}\right.
                                       \eqno{(1.4)}
$$
}
\\

Now we formulate \\
{\bf Backward problem:}\\
Let $T>0$ and $a_T, b_T$.  Then determine $u=u(x,t)$ such that 
$$
\left\{ \begin{array}{rl}
& \pppa u = -Au, \quad x\in \OOO, \,
t>0, \\
& u(\cdot,T) = a_T, \quad \ppp_tu(\cdot,T) = b_T, \qquad x \in \OOO, \\
& u(\cdot,t) \in H^1_0(\OOO), \quad t>0.
\end{array}\right.
$$

We set 
$$
\psi(\eta) := \MLO(-\eta)^2 + \eta \MLT(-\eta)\MLAA(-\eta), \quad
\eta > 0.                                    \eqno{(1.5)}
$$
By the definition of the Mittag-Leffler function, we have
$\psi(0) =  1$.
Before stating the main result, we show 
\\
{\bf Lemma 1.}\\
{\it 
The set $\{ \eta>0;\, \psi(\eta)=0\}$ is a non-empty and finite set.
}
\\

We set
$$
\{ \eta_1, ..., \eta_N\} = \{ \eta>0;\, \psi(\eta)=0\}
                                                  \eqno{(1.6)}
$$
with $\eta_1 < \cdots < \eta_N$.
We have no information of the number $N$ of the zeros of $\psi$, except 
that it exists.  In Lemma 2 in Section 4, we will provide an upper bound of
the largest zero $\eta_N$.

Now we are ready to state our main result:
\\
{\bf Theorem.}\\
{\it
(i) We assume 
$$
T \not\in \bigcup_{n=1}^{\infty} \left\{
\left( \frac{\eta_1}{\mu_n}\right)^{\frac{1}{\alpha}}, \cdots,
\left( \frac{\eta_N}{\mu_n}\right)^{\frac{1}{\alpha}} \right\}.
                                                   \eqno{(1.7)}
$$
Then for any $a_T, b_T \in H^2(\OOO) \cap H^1_0(\OOO)$, there exist
$a,b\in L^2(\OOO)$ such that the solution $u_{a,b}$ to 
(1.3) satisfies 
$$
\uab(\cdot,T) = a_T, \quad \ppp_t\uab(\cdot,T) = b_T.  \eqno{(1.8)}
$$
Moreover there exists a constant $C>0$ such that 
$$
C^{-1}(\Vert a_T\Vert_{H^2(\OOO)} + \Vert b_T\Vert_{H^2(\OOO)}) 
\le \Vert a\Vert_{L^2(\OOO)} + \Vert b\Vert_{L^2(\OOO)}
\le C(\Vert a_T\Vert_{H^2(\OOO)} + \Vert b_T\Vert_{H^2(\OOO)}) 
                                \eqno{(1.9)}
$$
for all $a_T, b_T \in L^2(\OOO)$.
\\
(ii) We assume
$$
T \in \bigcup_{n=1}^{\infty} \left\{
\left( \frac{\eta_1}{\mu_n}\right)^{\frac{1}{\alpha}}, \cdots,
\left( \frac{\eta_N}{\mu_n}\right)^{\frac{1}{\alpha}} \right\}.
                                          \eqno{(1.10)}
$$
Then there exists $(a, b) \not\equiv (0,0)$ in $\OOO$ such that 
$(u_{a,b}(\cdot,T), \, \ppp_tu_{a,b}(\cdot,T)) \equiv (0,0)$ in $\OOO$.
Furthermore, if $(u_{a,b}(\cdot,T), \, \ppp_tu_{a,b}(\cdot,T)) \equiv (0,0)$ 
in $\OOO$, then 
$$
(a, \va_{nj}) = (b, \va_{nj}) = 0, \quad 1\le j \le \ell_n
$$
if  
$$
T \not\in \left\{
\left( \frac{\eta_1}{\mu_n}\right)^{\frac{1}{\alpha}}, \cdots,
\left( \frac{\eta_N}{\mu_n}\right)^{\frac{1}{\alpha}} \right\}.
$$
}
\\

Henceforth we set 
$$
\Lambda = \Lambda(\alpha,A) 
:= \bigcup_{n=1}^{\infty} \left\{
\left( \frac{\eta_1}{\mu_n}\right)^{\frac{1}{\alpha}}, \cdots,
\left( \frac{\eta_N}{\mu_n}\right)^{\frac{1}{\alpha}} \right\}.
$$
We note that $\Lambda$ is a countably infinite set.
Theorem (i) implies that the backward problem in time for 
$1<\alpha<2$, is well-posed for some values of $T$ not belonging
to the non-empty set $\Lambda$.
The part (ii) means that we cannot determine the $\va_{nj}$-components of 
initial values where $n\in \N$ for which  
$T \in \left\{
\left( \frac{\eta_1}{\mu_n}\right)^{\frac{1}{\alpha}}, \cdots,
\left( \frac{\eta_N}{\mu_n}\right)^{\frac{1}{\alpha}} \right\}$, that is,
such exceptional values of the final time $T$ actually cause
the non-uniqueness for the backward problem.

Since $\lim_{n\to\infty} \mu_n = \infty$, the set $\Lambda$
has an accumulation point $0$, but we can readily verify
$$
\Lambda 
\subset \left[ 0, \left( \frac{\eta_N}{\mu_1}\right)^{\frac{1}{\alpha}}
\right].
$$
Hence
\\
{\bf Corollary 1.}\\
{\it 
If 
$$
T > \left( \frac{\eta_N}{\mu_1}\right)^{\frac{1}{\alpha}}, \eqno{(1.11)}
$$
then for any $a_T, b_T \in H^2(\OOO) \cap H^1_0(\OOO)$,
there exist unique $a,b,\in L^2(\OOO)$ such that 
$u_{a,b}$ satisfies (1.8) and (1.9).
}
\\

In Corollary 2 in Section 4, we provide a more concrete estimate of $T$ 
than (1.11). 

The backward problem for $1<\alpha<2$ is rather different from 
the case $0<\alpha<1$ which is well-posed for any $T>0$.
We can sum up the results for the backward problems for $0<\alpha \le 2$:
\\
{\bf Backward problem in time.}
\begin{itemize}
\item
{\bf $0<\alpha<1$:} well-posed for any $T>0$.
\item
{\bf $\alpha=1$:}
severely ill-posed but we have the uniqueness and some conditional stability
for any $T>0$.
\item
{\bf $1<\alpha<2$:}
well-posed for $T>0$ not belonging to a countably infinite set.
Even non-uniqueness occurs for such exceptional values of $T$.
\item
{\bf $\alpha=2$:}
Well-posed.  Also we have conservation quantity such as energy, which 
is impossible for $\alpha \ne 2$.
\end{itemize}
The well-posedness is sensitive according to $0<\alpha<1$, $\alpha=1$,
$1<\alpha<2$ and $\alpha=2$, and in the case $1<\alpha<2$, a quite new 
aspect of the non-uniqueness happens by choices of $T$.

This article is composed of four sections.  In Section 2, we prove 
Lemma 1, and Section 3 is devoted to the proof of Theorem.
Section 4 gives two concluding remarks.
\section{Proof of Lemma 1}

We recall that $\psi(\eta)$ is defined by (1.5). By the analyticity of the 
Mittag-Leffler function, we see that $\psi(\eta)$ is analytic in $\eta > 0$ 
and continuous in $[0, \infty)$.  Moreover by the asymptotics of the 
Mittag-Leffler functions (e.g., Theorem 1.4 (pp.33-34) in \cite{Po}), we 
see that 
$$
\left\{ \begin{array}{rl}
& \MLO(-\eta) = \frac{1}{\Gamma(1-\alpha)}\frac{1}{\eta} 
+ O\left(\frac{1}{\eta^2}\right), \quad
\MLT(-\eta) = \frac{1}{\Gamma(2-\alpha)}\frac{1}{\eta} 
+ O\left(\frac{1}{\eta^2}\right), \\
& \MLAA(-\eta) = \frac{-1}{\Gamma(-\alpha)}\frac{1}{\eta^2} 
+ O\left(\frac{1}{\eta^3}\right) \quad \mbox{as $\eta \to \infty$}.
\end{array}\right.
                                             \eqno{(2.1)}
$$
Therefore
\begin{align*}
& \psi(\eta) = \MLO(-\eta)^2 + \eta \MLT(-\eta)\MLAA(-\eta)\\
=& \left( \frac{1}{\Gamma(1-\alpha)}\frac{1}{\eta} 
+ O\left(\frac{1}{\eta^2}\right)\right)^2
- \eta\left( \frac{1}{\Gamma(2-\alpha)}\frac{1}{\eta} 
+ O\left(\frac{1}{\eta^2}\right)\right)
 \left(\frac{1}{\Gamma(-\alpha)}\frac{1}{\eta^2} 
+ O\left(\frac{1}{\eta^3}\right)\right)\\
=&  \left(\frac{1}{\Gamma(1-\alpha)^2} - \frac{1}{\Gamma(2-\alpha)
\Gamma(-\alpha)} \right) \frac{1}{\eta^2} 
+ O\left(\frac{1}{\eta^3}\right) \quad \mbox{as $\eta \to \infty$}.
\end{align*}
Since $\Gamma(1-\alpha) = -\alpha\Gamma(-\alpha)$ and 
$\Gamma(2-\alpha) = (1-\alpha)\Gamma(1-\alpha) = 
(\alpha^2-\alpha)\Gamma(-\alpha)$, we obtain
$$
\psi(\eta) = \frac{-1}{\alpha^2(\alpha-1)\Gamma(-\alpha)^2} \frac{1}{\eta^2} 
+ O\left(\frac{1}{\eta^3}\right) \quad \mbox{as $\eta \to \infty$.}
                                         \eqno{(2.2)}
$$
By $\frac{-1}{\alpha^2(\alpha-1)\Gamma(-\alpha)^2} < 0$, 
there exists a constant $M>0$ such that 
$\psi(\eta) < 0$ for $\eta \ge M$.
Since $\psi(0) = 1$, by the continuity of $\psi$ in $[0,\infty)$, we can 
find a sufficiently small constant 
$\ep>0$ such that $\psi(\eta) > 0$ for $0\le \eta \le \ep$.
Therefore the intermediate value theorem yields that there exists
$\eta_0 \in (\ep, M)$ such that $\psi(\eta_0) = 0$.
Moreover, since $\psi$ is analytic in $[\ep, M]$, the set 
$\{\eta \in [\ep,M];\, \psi(\eta) = 0\}$ is a finite set.  
Otherwise $\psi(\eta) = 0$ for 
each $\eta \in [\ep,M]$, which implies $\psi(0) = 0$ by the continuity of
$\psi(\eta)$ at $\eta = 0$, which contradicts $\psi(0) = 1$.
Thus the proof of Lemma 1 is complete.
\section{Proof of Theorem}

We set 
$$
a_{nj} = (a, \va_{nj}), \quad b_{nj} = (b, \va_{nj}), 
$$
and 
$$
\left\{ \begin{array}{rl}
& p_{nj}:= a_{nj}\MLO(-\mu_nT^{\alpha}) 
  + b_{nj}T\MLT(-\mu_nT^{\alpha}), \\
& q_{nj}:= -\mu_nT^{\alpha-1}a_{nj}\MLAA(-\mu_nT^{\alpha}) 
  + b_{nj}\MLO(-\mu_nT^{\alpha}).
\end{array}\right.
                                    \eqno{(3.1)}
$$
Since $\{\va_{nj}\}_{1\le j\le \ell_n, n\in \N}$ is an orthonormal 
basis in $L^2(\OOO)$, we see that 
\begin{align*}
&\sumn \suml \vert (g,\va_{nj})\vert^2 = \Vert g\Vert^2_{L^2(\OOO)},
\quad g \in L^2(\OOO), \\
&\sumn \suml \mu_n^2\vert (g,\va_{nj})\vert^2 = \Vert g\Vert^2_{H^2(\OOO)},
\quad g \in \mathcal{D}(A) = H^1(\OOO) \cap H^1_0(\OOO).
\end{align*}
Hence, by (1.4), we have
$$
\Vert u(\cdot,T)\Vert^2_{H^2(\OOO)} 
+ \Vert \ppp_tu(\cdot,T)\Vert^2_{H^2(\OOO)}     \eqno{(3.2)}
$$
\begin{align*}
=& \sumn \suml \mu_n^2(\vert a_{nj}\MLO(-\mu_nT^{\alpha}) 
  + b_{nj}T\MLT(-\mu_nT^{\alpha})\vert^2 \\
+& \vert -\mu_nT^{\alpha-1}a_{nj}\MLAA(-\mu_nT^{\alpha}) 
  + b_{nj}\MLO(-\mu_nT^{\alpha}) \vert^2)\\
=& \sumn \suml \mu_n^2(\vert p_{nj}\vert^2 + \vert q_{nj}\vert^2).
\end{align*}
Now we proceed to 
\\
{\bf Proof of Theorem (i).}\\
We assume (1.7), and so $\psi(\mu_nT^{\alpha}) \ne 0$ for all 
$n \in \N$.  Then we can solve (3.1) with respect to
$a_{nj}$ and $b_{nj}$:
$$
\left\{ \begin{array}{rl}
& a_{nj} = \frac{1}{\psi(\mu_nT^{\alpha})}
(p_{nj}\MLO(-\mu_nT^{\alpha}) - q_{nj}T\MLT(-\mu_nT^{\alpha})),\\
& b_{nj} = \frac{1}{\psi(\mu_nT^{\alpha})}
(p_{nj}\mu_nT^{\alpha-1}\MLAA(-\mu_nT^{\alpha}) 
                                 + q_{nj}\MLO(-\mu_nT^{\alpha})).
\end{array}\right.
                                    \eqno{(3.3)}
$$
By (2.1) and (2.2), we can choose a large constant $M_0>0$ such that 
\begin{align*}
& \vert \MLO(-\eta)\vert \le \frac{2}{\eta}
\left\vert \frac{1}{\Gamma(1-\alpha)}\right\vert,
\quad \vert \MLT(-\eta)\vert \le \frac{2}{\eta}\frac{1}{\Gamma(2-\alpha)},
\quad \vert \MLAA(-\eta)\vert \le \frac{2}{\eta^2}\frac{1}{\Gamma(-\alpha)},\\
& \vert \psi(\eta)\vert \ge \frac{1}{2\eta^2}\frac{1}{\alpha^2(\alpha-1)
\Gamma(-\alpha)^2}, \quad \eta \ge M_0.
\end{align*}
Here we note that $\Gamma(1-\alpha) < 0$ and $\Gamma(2-\alpha), 
\Gamma(-\alpha) > 0$.
Consequently we can fix $N_0 \in \N$ such that 
$$
\left\{ \begin{array}{rl}
& \vert \psi(\mu_nT^{\alpha})\vert 
\ge \frac{1}{2T^{2\alpha}\mu_n^2}\frac{1}{\alpha^2(\alpha-1)
\Gamma(-\alpha)^2} = \frac{C_1}{\mu_n^2}, \\
& \vert \MLO(-\mu_nT^{\alpha})\vert, \quad \vert \MLT(-\mu_nT^{\alpha})\vert,
\quad \vert \mu_n\MLAA(-\mu_nT^{\alpha})\vert \le \frac{C_1}{\mu_n}, 
\quad n \ge N_0.           
\end{array}\right.                         \eqno{(3.4)}
$$
Here and henceforth $C_k$, $k=1,2,..., 5,6$ denote generic constants 
which are independent of $n$ and $j$, but dependent on $T, N_0,
\alpha$.

Therefore (3.3) implies
$$
\vert a_{nj}\vert \le C_2\mu_n(\vert p_{nj}\vert + \vert q_{nj}\vert),
\quad 
\vert b_{nj}\vert \le C_2\mu_n(\vert p_{nj}\vert + \vert q_{nj}\vert),
\quad n\ge N_0, \, 1\le j \le \ell_n.      \eqno{(3.5)}
$$
On the other hand, by bounds of the
the Mittag-Leffler functions (e.g., Theorem 1.6 (p.35) in 
\cite{Po}), we see that 
$$
\vert \MLO(-\mu_nT^{\alpha})\vert, \quad 
\vert \MLT(-\mu_nT^{\alpha})\vert \le \frac{C_3}{1+\mu_n} \le C_4,
\quad n\in \N.                                        \eqno{(3.6)}
$$
Moreover the estimate of $\vert \MLAA(-\mu_nT^{\alpha})\vert$ in 
(3.4) implies
$$
\vert \MLAA(-\mu_nT^{\alpha})\vert \le \frac{C_4}{1+\mu_n^2},
\quad n\in \N.                             \eqno{(3.7)}
$$

Since $\psi(\mu_nT^{\alpha}) \ne 0$ for each $n \in \N$, by 
(3.3) and (3.6), we have
\begin{align*}
&\vert a_{nj}\vert \le C_5\max_{1\le n\le N_0-1}
\left\vert \frac{1}{\psi(\mu_nT^{\alpha})}\right\vert 
(\vert p_{nj}\vert + T\vert q_{nj}\vert), \\
&\vert b_{nj}\vert \le C_5\max_{1\le n\le N_0-1}
\left\vert \frac{1}{\psi(\mu_nT^{\alpha})}\right\vert 
(\mu_{n}T^{\alpha-1}\vert p_{nj}\vert + \vert q_{nj}\vert), 
\quad 1\le n\le N_0-1, \, 1\le j \le \ell_n,
\end{align*}
so that (3.5) holds for each $n \in \N$ and $1\le j \le \ell_n$.
Hence
$$
\sumn \suml (\vert a_{nj}\vert^2 + \vert b_{nj}\vert^2)
\le C_5\sumn \suml \mu_n^2(\vert p_{nj}\vert^2 + \vert q_{nj}\vert^2),
$$
and applying (3.2), we obtain
$$
\Vert a\Vert^2_{L^2(\OOO)} + \Vert b\Vert^2_{L^2(\OOO)}
\le C_5(\Vert u(\cdot,T)\Vert^2_{H^2(\OOO)} 
+ \Vert \ppp_tu(\cdot,T)\Vert^2_{H^2(\OOO)}).
$$
Next we prove the reverse inequality. 
Applying (3.6) and (3.7) in (3.1), we have
\begin{align*}
&\mu_n\vert p_{nj}\vert \le C_6( \vert a_{nj}\vert + \vert b_{nj}\vert),\\
&\mu_n\vert q_{nj}\vert \le C_6'\left( \frac{C_3}{1+\mu_n}\vert a_{nj}\vert
+ \vert b_{nj}\vert\right)
\le C_6(\vert a_{nj}\vert + \vert b_{nj}\vert)
\end{align*}
for all $n\in \N$ and $1\le j \le \ell_n$.  Hence, in view of (3.2), we see
$$
\Vert u(\cdot,T)\Vert^2_{H^2(\OOO)} + \Vert \ppp_tu(\cdot,T)\Vert^2
_{H^2(\OOO)}
\le C_6 \sumn \suml (\vert a_{nj}\vert^2+\vert b_{nj}\vert^2)
= C_6(\Vert a\Vert^2_{L^2(\OOO)} + \Vert b\Vert^2_{L^2(\OOO)}),
$$
which completes the proof of Theorem (i).
\\
{\bf Proof of Theorem (ii).}\\
By (1.4) we see
$$
u(\cdot,T) = \sumn \suml p_{nj}\va_{nj}, \quad
\ppp_tu(\cdot,T) = \sumn \suml q_{nj}\va_{nj}.
$$
Therefore $ u(\cdot,T) = \ppp_tu(\cdot,T) = 0$ in $\OOO$ is equivalent to
$p_{nj} = q_{nj} = 0$ for $n\in \N$ and $1\le j \le \ell_n$.
By (1.10), we can choose $n_0\in \N$ and $k_0 \in \{1, ..., N\}$ such that 
$T = \left(\frac{\eta_{k_0}}{\mu_{n_0}}\right)^{\frac{1}{\alpha}}$, that is,
$\eta_{k_0} = \mu_{n_0}T^{\alpha}$.  Consequently $\psi(\mu_{n_0}T^{\alpha}) 
= 0$.  Recalling the definition of $\psi(\mu_{n_0}T^{\alpha})$, we see that 
it is the determinant of the coefficient matrix of the linear system 
(3.1) with respect to $a_{n_0j}$ and $b_{n_0j}$.  Hence there exist
$(a_{n_01}, b_{n_01}) \ne (0,0)$ satisfying
$$
\left\{ \begin{array}{rl}
& a_{n_01}\MLO(-\mu_{n_0}T^{\alpha}) 
  + b_{n_01}T\MLT(-\mu_{n_0}T^{\alpha}) = 0, \\
& -\mu_{n_0}T^{\alpha-1}a_{n_01}\MLAA(-\mu_{n_0}T^{\alpha}) 
  + b_{n_01}\MLO(-\mu_{n_0}T^{\alpha}) = 0.
\end{array}\right.
$$
Setting $a= u(\cdot,0) := a_{n_01}\va_{n_01}$ and
$b= \ppp_tu(\cdot,0) := b_{n_01}\va_{n_01}$, we see that either 
$a \ne 0$ in $\OOO$ or $b\ne 0$ in $\OOO$, and
$u_{a,b}(\cdot,T) = \ppp_tu_{a,b}(\cdot,T) = 0$ in $\OOO$.  The former 
part of (ii) is now proved.  The latter part follows from (3.1).
Indeed let $T \not\in \left\{ \left(\frac{\eta_1}{\mu_n}\right)
^{\frac{1}{\alpha}}, ..., \left(\frac{\eta_N}{\mu_n}\right)
^{\frac{1}{\alpha}}\right\}$ for some $n\in \N$.  
Then $\psi(\mu_nT^{\alpha}) \ne 0$.
Therefore the determinant $\psi(\mu_nT^{\alpha})$ of the coefficient matrix 
of (3.1) is not zero, and so $a_{nj} = b_{nj} = 0$, that is,
$(a,\va_{nj}) = (b, \va_{nj}) = 0$ for $1\le j\le \ell_n$.
Thus the proof of Theorem is complete.  

\section{Concluding remarks}

{\bf 4.1. Estimation of $\eta_N$ and $T$.}

We recall (1.6).  First we give an upper bound for $\eta_N$.
For simplicity, we set 
\begin{align*}
&\mu_1 := \frac{-1}{\Gamma(1-\alpha)} = \frac{1}{\alpha\Gamma(-\alpha)},\\
&\mu_2 := \frac{1}{\Gamma(2-\alpha)} 
= \frac{1}{(\alpha^2-\alpha)\Gamma(-\alpha)}, \quad 
\mu_3:= \frac{1}{\Gamma(-\alpha)}.
\end{align*}
Here we used $\Gamma(1-\alpha) = -\alpha\Gamma(-\alpha)$ and
$\Gamma(2-\alpha) = (1-\alpha)\Gamma(1-\alpha) 
= (1-\alpha)(-\alpha)\Gamma(-\alpha)$.
By $\Gamma(-\alpha) > 0$, we see that $\mu_1, \mu_2, \mu_3 > 0$.

For $1<\alpha<2$, we choose $\theta$ such that 
$$
\frac{\pi\alpha}{2} < \theta < \pi.            \eqno{(4.1)}
$$
By $\gamma$ we denote the contour in $\C$ which is directed from 
$\infty e^{-\sqrt{-1}\theta}$ to $\infty e^{\sqrt{-1}\theta}$ and 
consists of \\
(i) arg $z = -\theta$, $\vert z\vert \ge 1$ \\
(ii) $-\theta \le \mbox{arg}\, z \le \theta$, $\vert z\vert = 1$ \\
(iii) arg $z = \theta$, $\vert z\vert \ge 1$. 

Moreover we set 
\begin{align*}
& \nu_1 = \frac{1}{2\pi \alpha \sin\theta}\int_{\gamma}
\vert \exp(\zeta^{\frac{1}{\alpha}}) \vert \vert \zeta\vert d\zeta,\\
& \nu_2 = \frac{1}{2\pi \alpha \sin\theta}\int_{\gamma}
\vert \exp(\zeta^{\frac{1}{\alpha}}) \vert \vert \zeta^{1-\frac{1}{\alpha}}
\vert d\zeta,\\
& \nu_3 = \frac{1}{2\pi \alpha \sin\theta}\int_{\gamma}
\vert \exp(\zeta^{\frac{1}{\alpha}}) \vert \vert \zeta^{1+\frac{1}{\alpha}}
\vert d\zeta.
\end{align*}
Since there exists a constant $C_0>0$ such that 
$\eexx \le \exp\left(-C_0\vert \zeta\vert^{\frac{1}{\alpha}}\right)$
for $\zeta \in \gamma$, we can directly verify that 
$0 < \nu_1, \nu_2, \nu_3 < \infty$.

Then we can prove
\\
{\bf Lemma 2.}\\
{\it
$$
\eta_N <  \max \biggl\{ \frac{1}{\vert \cos \theta\vert}, 
\, \alpha^2(\alpha-1)\Gamma(-\alpha)^2
(\mu_2\nu_3 + \mu_3\nu_2 + 2\mu_1\nu_1 + \nu_1^2 + \nu_2\nu_3)
\biggr\}.
$$
}
\\
{\bf Proof of Lemma 2.}\\
First by formula (1.145) (p.34) in \cite{Po}, we see that 
$$\left\{ \begin{array}{rl}
& \MLO(-\eta) = - \frac{\mu_1}{\eta} +I_{\alpha,1}(\eta), \quad
\MLT(-\eta) = \frac{\mu_2}{\eta} +I_{\alpha,2}(\eta), \\
& \MLAA(-\eta) = - \frac{\mu_3}{\eta^2} +I_{\alpha,\alpha}(\eta), \quad
\eta \ge 1,
\end{array}\right.                   \eqno{(4.2)}
$$
where 
$$
I_{\alpha,\ell}(\eta) = \frac{-1}{2\pi\alpha \sqrt{-1}\eta}
\int_{\gamma} \exp(\zeta^{\frac{1}{\alpha}}) \zeta^{\frac{1-\ell}{\alpha}+1}
\frac{d\zeta}{\zeta+\eta},
$$
$$ 
I_{\alpha,\alpha}(\eta) = \frac{1}{2\pi\alpha \sqrt{-1}\eta^2}
\int_{\gamma} \exp(\zeta^{\frac{1}{\alpha}}) \zeta^{\frac{1}{\alpha}+1}
\frac{d\zeta}{\zeta+\eta}, \quad \ell=1,2, \quad \eta \ge 1.       \eqno{(4.3)}
$$
Next we will prove
$$
\vert I_{\alpha,1}(\eta) \vert \le \frac{\nu_1}{\eta^2}, \quad
\vert I_{\alpha,2}(\eta) \vert \le \frac{\nu_2}{\eta^2}, \quad
\vert I_{\alpha,\alpha}(\eta) \vert \le \frac{\nu_3}{\eta^3} \quad
\mbox{for $\eta \ge \frac{1}{\vert \cos\theta\vert}$}.    \eqno{(4.4)}
$$
{\bf Proof of (4.4).}
For $\eta \ge \frac{1}{\vert \cos\theta\vert}$, we can 
directly verify that 
$$
\min_{\zeta\in \gamma} \vert \zeta + \eta\vert 
= \min_{r\ge 1} \vert re^{\sqrt{-1}\theta} - (-\eta)\vert 
= \vert -\eta\vert \sin (\pi - \theta) = \eta \sin\theta>0. 
$$
Indeed the intersection point of the perpendicular from $-\eta$ with 
the half-line $\{ re^{\sqrt{-1}\theta};\, r\ge 1\}$ is outside of 
$\{z\in \C;\, \vert z\vert \le 1\}$ if $\vert \eta\vert > \frac{1}
{\vert \cos\theta\vert}$.  Hence with fixed $\eta$ satisfying 
$\eta \ge \frac{1}{\vert \cos\theta\vert}$, 
the function $\vert \zeta+\eta\vert$ in $\zeta \in \gamma$, attains the minimum
at such an intersection point $\zeta$.

Therefore
\begin{align*}
& \vert I_{\alpha,\ell}(\eta) \vert 
\le \frac{1}{2\pi\alpha\sin\theta}\frac{1}{\eta^2}
\int_{\gamma} \eexx \vert \zeta^{\frac{1-\ell}{\alpha}+1}\vert d\zeta
= \nu_{\ell}\frac{1}{\eta^2}, \\
& \vert I_{\alpha,\alpha}(\eta) \vert 
\le \frac{1}{2\pi\alpha \sin\theta \eta^3}
\int_{\gamma} \eexx \vert \zeta^{\frac{1}{\alpha}+1}\vert d\zeta
= \frac{\nu_3}{\eta^3}, \quad
\eta \ge \frac{1}{\vert \cos\theta \vert}.
\end{align*}
Hence (4.4) is proved.
\\

Now we will complete the proof of Lemma 2.
Applying (4.2) and (4.4) in (1.5), we obtain
\begin{align*}
& \psi(\eta) = \MLO(-\eta)^2 + \eta\MLT(-\eta)\MLAA(-\eta)\\
=& \left( \frac{\mu_1}{\eta} - I_{\alpha,1}(\eta)\right)^2
+ \left( \frac{\mu_2}{\eta} +I_{\alpha,2}(\eta)\right)
\left(- \frac{\mu_3}{\eta} + \eta I_{\alpha,\alpha}(\eta)\right)\\
=& \frac{1}{\eta^2}(\mu_1^2 - \mu_2\mu_3)
+ \left\{ -2\frac{\mu_1}{\eta}I_{\alpha,1}(\eta)
+ I_{\alpha,1}(\eta)^2 + \mu_2I_{\alpha,\alpha}(\eta)
- \frac{\mu_3I_{\alpha,2}(\eta)}{\eta}
+ \eta I_{\alpha,2}(\eta)I_{\alpha,\alpha}(\eta)\right\}\\
\le& \frac{1}{\eta^2}(\mu_1^2 - \mu_2\mu_3)
+ \left( 2\frac{\mu_1}{\eta}\vert I_{\alpha,1}(\eta)\vert 
+ \vert I_{\alpha,1}(\eta)\vert^2 + \mu_2\vert I_{\alpha,\alpha}(\eta)\vert
+ \frac{\mu_3\vert I_{\alpha,2}(\eta)\vert}{\eta}
+ \eta \vert I_{\alpha,2}(\eta)\vert 
\vert I_{\alpha,\alpha}(\eta)\vert \right)\\
\le& \frac{-1}{\alpha^2(\alpha-1)\Gamma(-\alpha)^2}\frac{1}{\eta^2}
+ (\mu_2\nu_3 + \mu_3\nu_2 + 2\mu_1\nu_1)\frac{1}{\eta^3}
+ (\nu_1^2+\nu_2\nu_3)\frac{1}{\eta^4}\\
\le& \frac{-1}{\alpha^2(\alpha-1)\Gamma(-\alpha)^2}\frac{1}{\eta^2}
+ (\mu_2\nu_3 + \mu_3\nu_2 + 2\mu_1\nu_1 + \nu_1^2+\nu_2\nu_3)
\frac{1}{\eta^3}
\end{align*}
for $\eta \ge \frac{1}{\vert \cos\theta\vert} \ge 1$.  For the last
inequality, we use $\frac{1}{\eta^4} \le \frac{1}{\eta^3}$ for 
$\eta \ge 1$.  Hence
\begin{align*}
& \psi(\eta) \\
= & \frac{-1}{\alpha^2(\alpha-1)\Gamma(-\alpha)^2}\frac{1}{\eta^2}
\left\{ 1 - \alpha^2(\alpha-1)\Gamma(-\alpha)^2
(\mu_2\nu_3 + \mu_3\nu_2 + 2\mu_1\nu_1 + \nu_1^2+\nu_2\nu_3)\frac{1}{\eta}
\right\}
\end{align*}
for $\eta \ge \frac{1}{\vert \cos\theta\vert}$.  Thus the proof of Lemma 2 is
complete.
\\

Applying Lemma 2 to Corollary 1, we obtain a lower bound of $T$ which 
is described more concretely than (1.11) for guaranteeing 
the well-posedness f the backward problem.
\\
{\bf Corollary 2.}\\
{\it
If 
$$
T > \left( \frac{1}{\mu_1} \max \biggl\{ \frac{1}{\vert \cos \theta\vert}, 
\alpha^2(\alpha-1)\Gamma(-\alpha)^2
(\mu_2\nu_3 + \mu_3\nu_2 + 2\mu_1\nu_1 + \nu_1^2 + \nu_2\nu_3)
\biggr\}\right)^{\frac{1}{\alpha}},
$$
then for any $a_T, b_T \in H^2(\OOO) \cap H^1_0(\OOO)$,
there exist unique $a,b\in L^2(\OOO)$ such that 
$u_{a,b}$ satisfies (1.8) and (1.9).
}
\\

{\bf 4.2. Backward fractional ordinary differential equations.}

Let $1<\alpha<2$ and $\la > 0$.  We consider a backward fractional ordinary 
differential equation.
$$
\pppa v(t) = -\la v(t), \quad v(T) = a_T, \quad \ppp_tv(T)=b_T, \quad 0<t<T.
                                            \eqno{(4.5)}
$$
By (1.6), we can prove that $T \not\in \left\{
\left(\frac{\eta_1}{\la}\right)^{\frac{1}{\alpha}}, \, ...,
\left(\frac{\eta_N}{\la}\right)^{\frac{1}{\alpha}} \right\}$, then 
(4.5) possesses a unique solution for arbitrary $a_T, b_T \in \R$.
Moreover if $T \in \left\{
\left(\frac{\eta_1}{\la}\right)^{\frac{1}{\alpha}}, \, ...,
\left(\frac{\eta_N}{\la}\right)^{\frac{1}{\alpha}} \right\}$, then 
there exists a non-zero solution $v$ to (4.5) with 
$a_T=b_T=0$, and there may be no solutions with some $a_T, b_T$.
Thus for the case $1<\alpha<2$,
the backward problem for a fractional ordinary differential equation is
not always uniquely solvable for all $T>0$.  In general, even for 
nonlinear fractional ordinary differential equations, 
under suitable conditions, we can apply the contraction mapping theorem
to prove the well-posedness for sufficiently small $T$.  
However, as Theorem (ii) asserts, the backward problem for fractional 
partial differential equations with $1<\alpha<2$ may not be well-posed
even for sufficiently small $T>0$.  We can refer to an example 
(p.374) in Diethelm and Ford \cite{DF} which indicates the non-uniqueness 
in the case of $b_T=0$ with some value of $\la$.

\section*{Acknowledgment}
This work is also supported by the Istituto Nazionale di Alta Matematica 
(IN$\delta$AM),
through the GNAMPA Research Project 2019.  


The second author was supported by Grant-in-Aid for Scientific Research (S)
15H05740 of Japan Society for the Promotion of Science and
by The National Natural Science Foundation of China
(no. 11771270, 91730303).
This work was prepared with the support of the "RUDN University Program 5-100".


\end{document}